\documentclass{amsart}
\usepackage{amsfonts}
\usepackage{amssymb,amsmath,amscd,cancel}
\usepackage{epsfig}
\usepackage{graphicx}
\usepackage[all]{xy}

\newtheorem{thm}{\sc Theorem}[section]
\newtheorem{prop}[thm]{\sc Proposition}

\newtheorem{cor}[thm]{\sc Corollary}
\theoremstyle{definition}

\theoremstyle{definition}
\newtheorem{de}[thm]{\sc Definition}
\theoremstyle{definition}
\newtheorem{rem}[thm]{\sc Remark}
\theoremstyle{definition}

\numberwithin{equation}{section}

\begin{document}
\title[The fundamental theorem of curves and classifications]{The fundamental theorem of curves and classifications in the Heisenberg groups}

\author[H.-L.~Chiu, X.-H. Feng, Y.-C.Huang, ]{Hung-Lin Chiu, XiuHong Feng ,Yen-Chang Huang}
\address{Department of Mathematics, National Central University, Chung Li,
32054, Taiwan, R.O.C.}
\email{hlchiu@math.ncu.edu.tw}
\address{Department of Mathematics,The School of Mathematics and Statistics, Nanjing University of Information Science and Technology}
\email{fengxiuhong@nuist.edu.cn}
\address{Department of Mathematics, National Central University, Chung Li,
32054, Taiwan, R.O.C.}
\email{yenchang.huang1@gmail.com}

\begin{abstract}
We study the horizontally regular curves in the Heisenberg groups $H_n$. We show the fundamental theorem of curves in $H_n$ $(n\geq 2)$ and define the orders of horizontally regular curves. We also show that the curve $\gamma$ is of order $k$ if and only if, up to a Heisenberg rigid motion, $\gamma$ lies in $H_k$ but not in $H_{k-1}$; moreover, two curves with the same order differ from a rigid motion if and only if they have the same invariants: p-curvatures and contact normality. Thus, combining with our previous work \cite{CHL} we have completed the classification of horizontally regular curves in $H_n$ for $n\geq 1$.
\end{abstract}

\maketitle

\section{Introduction}
The Heisenberg group $H_n$, $n\geq 1$, is the space $\mathbb{R}^{2n+1}$ associated with the group multiplication
\begin{align*}
&(x_1,\cdots,x_n,y_1,\cdots, y_n,z)\circ(\widetilde{x_1},\cdots,\widetilde{x_n},\widetilde{y_1},\cdots,\widetilde{y_n},\widetilde{z})\\
=&(x_1+\widetilde{x_1},\cdots,x_n+\widetilde{x_n},y_1+\widetilde{y_1},\cdots,{y_n}+\widetilde{y_n},
z+\widetilde{z}+\sum_{j=1}^{n}(y_j\widetilde{x_j}-x_j\widetilde{y_j}))
\end{align*}
It is a $(2n+1)$-dimensional Lie group, and the space of all left invariant vector fields is spanned by the basic
vector fields:
\begin{align*}
\mathring{e}_j&= \frac{\partial}{\partial x_j}+{y_j}\frac{\partial}{\partial z}, \\
\mathring{e}_{n+j}&=\frac{\partial}{\partial y_j}-x_j\frac{\partial}{\partial z}, \\
T&=\frac{\partial}{\partial z},
\end{align*}
for $1\leq j\leq n.$

The Heisenberg group $H_n$ can be regarded as the $n$-dimensional CR manifold with zero Webster-curvature. For more details, the reader can refer the Appendix in \cite{CHMY} or \cite{CHL}\cite{CL}\cite{CHMY2}\cite{Jaco}. We give the brief description of geometric structures on $H_n$: the standard contact bundle in $H_n$ is the subbundle $\xi$ of the tangent bundle $TH_n$, which is spanned by $\mathring{e}_j$ and
$\mathring{e}_{n+j}$ for $1\leq j\leq n$. The contact bundle can also be defined as the kernel of the contact form $$\theta=dz+\sum_{j=1}^n (x_jdy_j-y_jdx_j).$$
The standard CR structure on $H_n$ is the almost complex structure defined on $\xi$ by
$$J(\mathring{e}_j)=\mathring{e}_{n+j},\ J\mathring{e}_{n+j}=-\mathring{e}_j.$$
Throughout the article, we regard the Heisenberg group $H_n$ with the standard pseudo-hermitian structure $(J,\theta)$
as a pseudohermitian manifold $(H_n, J,\theta)$. Denote the group of pseudohermitian transformations on $H_n$ by $PSH(n)$ which forms the group of rigid motions. An element in $PSH(n)$ is called a \textit{pseudohermitian transformation} or a \textit{symmetry} on $H_n$, which is a deffeomorphism $\Phi:H_n\rightarrow H_n$ preserving both the CR structure $J$ and the contact
form $\theta$. More precisely, it satisfies
$$\Phi_*J=J\Phi_*,\ \ \ \ \ \Phi^*\theta=\theta.$$

Here are our settings for curves: suppose $\gamma: I\subset \mathbb{R}\rightarrow H_n$ is a parametrized curve defined by
\begin{align}
\gamma(t) = (x_1(t),\cdots,x_n(t),y_1(t),\cdots,y_n(t),z(t)).
\end{align}
For $k=1,\cdots,n$, the $k$th derivative $\gamma^{(k)}$ of the curve $\gamma$ has the natural decomposition
\begin{align}\label{decomp}
\gamma^{(k)}(t) = \gamma^{(k)}_\xi (t)+ \gamma^{(k)}_T(t),
\end{align}
where $\gamma^{(k)}_\xi$ (resp. $\gamma^{(k)}_T$) is the orthogonal projection of $\gamma^{(k)}$ on the contact plane $\xi$ along $T$-direction (resp. on $T$ along $\xi$) with respect to the Levi metric. Recall that a curve is called \textit{horizontally regular} if it has the non-vanishing first derivative in the horizontal part, $\gamma'_\xi(t)\neq 0$ for all $t$. In \cite{CHL} we show that any horizontally regular curve can be reparametrized by horizontal arc length $s$ with respect to the Levi metric, namely, $|\gamma'_\xi(s)|=1$ for all $s$. In the article, we use different parameters $t$ and $s$ to distinguish from being parametrized by arc length or not. Moreover, by identifying $H_n\cong \mathbb{C}^n \times \mathbb{R}$ and the natural projection,
\begin{displaymath}
\xymatrixcolsep{4pc}\xymatrix{ & H_n\cong \mathbb{C}^n\times \mathbb{R}  \ar[d]^\pi  \\  I \ar[ur]^\gamma \ar[r]_\beta & \hspace{0.7cm}\mathbb{C}^n\cong \mathbb{R}^{2n}}
\end{displaymath}
we may rewrite the curve $\gamma$ in the real sense 
\begin{align}\label{curvehn}
\gamma(t)=\big(x_1(t),\cdots,x_n(t),y_1(t),\cdots, y_n(t),z(t)\big)\in \mathbb{R}^{2n}\times \mathbb{R},
\end{align}
with its projection on $\mathbb{R}^{2n}$
\begin{align}\label{rproj}
\alpha(t)=(x_1(t),y_1(t),\cdots, x_n(t),y_n(t)),
\end{align}
or equivalently in the complex sense
\begin{align}\label{cproj}
\beta(t)=\big(z_1(t),\cdots,z_n(t)\big)\in \mathbb{C}^n,
\end{align}
where
$z_j(s)=x_j(s)+\sqrt{-1}y_j(s)$ for $1\leq j\leq n$.

A key observation is that in case $H_1\ni \gamma=(\beta, z)$, we have
\begin{align*}
\gamma'& = (\beta', z') \in \mathbb{C}\times \mathbb{R}  \\
&=(x',y',z') \in \mathbb{R}^2\times \mathbb{R}\\
&=x'\mathring{e}_1 + y'\mathring{e}_2 +(z'-yx'+x'y)T.
\end{align*}
Thus, $\beta'\neq 0$ if and only if the curve $\gamma$ is horizontally regular in $H_1$. Given a curve in $H_n$, in general, we ask if one can establish the concept such that the curve is horizontally regular in \textit{any lower dimensional} subspaces of $H_n$.

Recall \cite{Jaco} a real linear subspace $P$ is totally real if and only if any vector $X\in P$ implies $JX\notin P$. Inspired by Griffiths \cite{Gri}, we generalize the definition of non-degenerate curves in $H_n$.
\begin{de}\label{nondeg}
\textbf{A non-degenerate horizontally regular curve} in $H_n$ is a horizontally regular curve $\gamma(t)=(\beta(t),z(t))\in\mathbb{C}^n\times \mathbb{R}$ satisfying
\begin{align}\label{cond1}
W_\beta^{[n]}(t):=&\gamma_\xi'(t)\wedge\cdots\wedge \gamma_\xi^{(n)}(t)\neq 0 \text{ for all } t,
\end{align}
and the set
\begin{align}\label{cond2}
\{ \gamma_\xi'(t),\cdots,\gamma_\xi^{(n)}(t)\} \text{ is totally real for all }t.
\end{align}
%A point $\gamma(t_0)\in H_n$ is called horizontally regular if $W_\beta^{[1]}(t_0)\neq 0$. Actually the curve $\gamma(t)\in H_n$ is horizontally regular if $W^{[1]}_\beta(t)\neq 0 \mbox{ for all } t$.
\end{de}

The condition (\ref{cond1}) ensures that we can always choose an oriented frame along the non-degenerate horizontally regular curve $\gamma$
\begin{align}\label{frame}
(\gamma(s); e_1(s),\cdots ,e_n(s),e_{n+1}(s),\cdots,e_{2n}(s),T)
\end{align}
inductively satisfying the condition, for $1\leq j\leq n$,
$$e_1(s)\wedge \cdots\wedge e_j(s)=\pm \frac{\gamma_{\xi}^{'}(s)\wedge \cdots \wedge
\gamma_{\xi}^{(j)}(s)}{|\gamma_{\xi}^{'}(s)\wedge \cdots \wedge \gamma_{\xi}^{(j)}(s)|}$$ and
$$ e_{n+j}=Je_j.$$ Along the curve $\gamma$, as in \cite{CHL}, we define the \textbf{p-curvatures}
$\kappa_j(s)$, $1\leq j\leq n$ and the \textbf{contact normality} $ \tau(s)$ as
\begin{align*}
\begin{array}{ll}
\displaystyle{\kappa_j(s)=\langle\frac{de_j(s)}{ds},e_{n+j}(s)\rangle, \text{ for } 1\leq j\leq n-1}, \\
\displaystyle{\kappa_n(s)=\langle\frac{de_n(s)}{ds},e_{2n}(s)\rangle},\\
\displaystyle{\tau(s)=\langle\frac{d\gamma(s)}{ds},T\rangle}.
\end{array}
\end{align*}
We point out that all quantities above are invariant under the group actions of $PSH(n)$. Our main theorem shows that those invariants completely determine the non-degenerate horizontally regular curve up to a Heisenberg rigid motion, which is analogues to the fundamental theorem of curves in $\mathbb{R}^n$.

\begin{thm}\label{mthm1}
Given $(n+1)$ smooth functions $\kappa_i(s)$ for $1\leq i\leq n$ and $\tau(s)$, there exists a non-degenerate horizontally regular curve $\gamma(s)\in H_n$ parametrized by the horizontal arc-length $s$ such that the functions $\kappa_i$'s and
$\tau$ are the p-curvatures and the contact normality of $\gamma$, respectively.
In addition, two non-degenerate horizontally regular curves satisfying the same conditions above differ from a rigid motion in $PSH(n)$.
\end{thm}

It is obvious that $\gamma(s)$ is a horizontal if $\gamma^{'}(s)=\gamma^{'}_{\xi}(s)$ for all $t \in I$, and therefore $\gamma(s)$ is horizontal if and only if $\tau(s)=0$. We immediately have the corollary:

\begin{cor}
Given smooth functions $\kappa_i(s),1\leq i\leq n,$ there exists a
horizontal curve $\gamma(s)\in H_n$ parametrized by the horizontal arc-length $s$ having $\kappa_i(s),1\leq i\leq n$  as its p-curvatures.
In addition, two horizontal curves having the same p-curvatures differ from a rigid motion in $PSH(n)$.
\end{cor}

Next we define the order of the curves. It is similar to the concept that spacial curves in $\mathbb{R}^3$ can not be "squeezed" into any linear 2-dimensional subspaces, but planar curves can be. The \textit{order} of a horizontal curve gives the minimal dimension of the subspaces in which the curve lives.

\begin{de}\label{deg}
A horizontally regular curve $\gamma(t)=(\beta(t),z(t))\in H_n$ is of \textbf{order k}, denoted by \textit{order}$(\gamma)=k$, if there exists a positive integer $k\in [1,n]$ such that
\begin{align}\label{order}
\left\{
\begin{array}{ll}
&\beta'(t)\wedge \cdots \wedge \beta^{(k+1)}(t)=0, \\
&\beta'(t)\wedge \cdots \wedge \beta^{(k)}(t)\neq 0,
\end{array}\right.
\end{align} for all $t$.
A horizontally regular curve is called \textbf{degenerate} in $H_n$ if \text{order}$(\gamma)< n$.
\end{de}

\begin{rem}
By Definition \ref{nondeg}, any non-degenerate horizontally regular curve $\gamma$ is of the top order, \textit{order}$(\gamma)=n$, and vice versa. In contrast to Theorem \ref{mthm1}, two curves with different orders never lie in \emph{the same subspace} of $H_n$, and hence they can not be congruent to each other by any Heisenberg rigid motion.
\end{rem}

We also characterize the degenerate horizontally regular curves of top order $(n-1)$. Similar to the fact that a planar curve in $\mathbb{R}^3$ can be "moved" to $xy$-plane, any degenerate regular curve $\gamma\in H_n$ can be acted by a symmetry of $PSH(n)$ to $H_{n-1}$.

\begin{prop}\label{topdeg}
Let $\gamma(t)=(\beta(t),z(t))\in H_n\cong \mathbb{C}^n\times \mathbb{R}$ be a degenerate horizontally regular curve. Then there exists a symmetry $\Phi\in PSH(n)$ such that $\Phi(\gamma)=\widetilde{\gamma}$, where $\widetilde{\gamma}=(\widetilde{\beta}(t),z(t))$ is a horizontally regular curve with the projection $\widetilde{\beta}=(\widetilde{\beta_1},\cdots, \widetilde{\beta_{n-1}},0)$ of $\beta$ onto $\mathbb{C}^{n-1}$. Thus, we conclude that $\widetilde{\gamma} \in H_{n-1}\subset H_n $.
\end{prop}

\begin{rem}\label{summary}
In summary, we have the dichotomy to classify any horizontally regular curve $\gamma$ in $H_n$: let $\gamma(t)=(\beta(t),z(t))\in H_n$. If $\gamma$ is non-degenerate, then by Theorem \ref{mthm1} it must be uniquely determined by the $p$-curvatures $\kappa_i(s)$, $1\leq i\leq n$, and contact normality $\tau(t)$ up to a symmetry in $PSH(n)$; otherwise, the Wronskian $W_\beta^{[n]}(t)= 0$ somewhere and we keep checking if $W_\beta^{[n-1]}(t)$ is nonzero for all $t$. The nonzero condition implies that $\gamma$ is degenerate of order $(n-1)$. By using Proposition \ref{topdeg} and applying Theorem \ref{mthm1} to $H_{n-1}$, we obtain that $\gamma$ lies in the subspace $H_{n-1}$ but not in $H_{n-2}$, and is uniquely determined by the invariants $\kappa_1,\cdots,\kappa_{n-1}, \tau$. However, if $W_\beta^{[n-1]}(t)=0$ somewhere, we have to check if $W_\beta^{[n-2]}(t)$ is nonzero again for all $t$. Repeating above processes, and finally we conclude that two curves with the same order $k\leq n-1$ differ from a symmetry in $PSH(k)$ if and only if both have same $k_i$, $i=1,\cdots, k-1$ and $\tau$. Thus, we complete the classification of horizontally regular curves.
\end{rem}

An interesting example is that the order of horizontal geodesics in $H_n$, $n\geq 1$, is always $1$. By the processes described in Remark \ref{summary}, the horizontal geodesics can always be embedded into $H_1$.

\begin{prop}\label{geodesic}
Every horizontal geodesic $\gamma\in H_n$, $n\geq 1$, is of order $1$ with constant $p$-curvatures and zero contact normality. Therefore, $\gamma$ is non-degenerate in $H_1$ and degenerate in $H_n$ for $n\geq 2$.
\end{prop}

The structure of the paper is as follows: in Section 2 we recall the well-known theorems for existence and uniqueness. In Section 3 we derive the Darboux derivatives. In Section 4 we prove the Theorem \ref{mthm1}. Finally, in Section 5, we characterize the degenerate curves (Proposition \ref{topdeg}), and, as an example, the order of horizontal geodesics will be calculated.

\textbf{Acknowledge} The first author would like to thank the Ministry of Science and Technology for the support of the project: MOST-104-2115-M-008-003-MY2, and also thanks for support from NCTS. The third author would like to express his appreciation to Professor Paul Yang and Professor Sun-Yung Alice Chang for their comments and suggestions.

\section{Calculus on Lie groups}
In the section, we shall give two well-known and essentially local results concerning smooth maps from manifold $M$ into the Lie group $G$. Given a connected smooth manifold $M$. Let $G \subset GL(n,R)$ be a matrix Lie group with Lie algebra $\mathfrak{g}$ and
the (left-invariant) Maurer-Cartan form $\omega$ on G. The first result is the existence theorem:

\begin{thm}[\cite{Gri}]\label{existence}
Suppose that $\phi$ is a $\mathfrak{g}$-valued one form on a simply connected manifold $M$.
Then there exists a $C^\infty$-map $f:M\rightarrow G$ with $f^*\omega=\phi$
if and only if $d\phi+\phi\wedge \phi=0$.
\end{thm}

The second result states that the pull-back of the Maurer-Cartan form uniquely determines the map up to a group action:
\begin{thm}[\cite{Gri}]\label{uniqueness}
Given two maps $f,\widetilde{f}:M\rightarrow G$, then $\widetilde{f}^*\omega=f^*\omega$ if and only if
$\widetilde{f}=g\circ f$ for some $g \in G$.
\end{thm}

We call the Lie algebra valued one-form $f^*\omega$ the \textit{Darboux derivative} of the map $f:M\rightarrow G$.

\section{ Differential Invariants of Horizontally Regular Curves in $H_n$}
Recall (equations (4.9)(4.10) in \cite{CL}) that any point $p\in H_n$ and element $\Phi_{p,R}\in PSH(n)$ have the corresponding representations respectively
\begin{align*}
p\in H_n &\longleftrightarrow X:=\begin{pmatrix}1 \\ p \end{pmatrix} \in\mathbb{R}^{2n+2}, \\
\Phi_{p,R}\in PSH(n) &\longleftrightarrow M\in PSH(n),
\end{align*}
satisfying the matrix multiplication
\begin{align*}
MX=\begin{pmatrix}1\\ \Phi_{p,R}\begin{pmatrix}p\end{pmatrix}\end{pmatrix}.
\end{align*}
Denote the indices
\begin{align*}
1\leq a&,\ b\leq 2n,\\
1\leq&\ j,k  \leq n.
\end{align*}
We also have the Maurer-Cartan form $\omega$ of $PSH(n)$ (\cite{CL}, page 1104)
\begin{align*}
\left(\begin{array}{cccc}0&0&0&0 \\ \omega^{k}&\omega_{j}^{k}&\omega_{n+j}^{k}&0
\\ \omega^{n+k}&\omega_{j}^{n+k}&\omega_{n+j}^{n+k}&0
\\ \omega^{2n+1}&\omega^{n+j}&-\omega^{j}&0
\end{array}\right),
\end{align*}
where
$\omega^{k},\ \omega_{j}^{k},\ \omega^{2n+1}$ are 1-forms on $PSH(n)$ satisfying
$\omega_a^b=-\omega^a_b$, $\omega_{j}^{n+k}=-\omega_{n+j}^{k}$, $\omega_{j}^{k}=\omega_{n+j}^{n+k}$.

Let $(p; e_j, e_{n+j}, T)$ be an oriented frame at point $p\in H_n$. By identifying $PSH(n)$ with the space of all oriented frames on $H_n$,
$$M \in PSH(n) \leftrightarrow (p,e_{j},e_{n+j},T),$$ we have
$$(p,e_{j},e_{n+j},T)=(0,\mathring{e}_{j},\mathring{e}_{n+j},T)\cdot M, $$ where $\cdot$ denotes the matrix multiplication. Thus, one can derive the moving frame formulas (page 1105, \cite{CL})
\begin{align}\label{mformula}
\begin{array}{rccccl}
dp  =&e_j\omega^j   &+&e_{n+j}\omega^{n+j}&+&T\omega^{2n+1}, \\
de_j=&e_k\omega^k_j &+&e_{n+k}\omega^{n+k}_j&+&T\omega^{n+j}, \\
de_{n+j}=&e_k\omega^k_{n+j}&+&e_{n+k}\omega^{n+k}_{n+j}&-&T\omega^j,\\
dT=&0. &&&&
 \end{array}
\end{align}

Let $\gamma(s)$ be a horizontally regular curve with horizontal arc-length parameter $s$. Each point of $\gamma$ uniquely defines the oriented frame as (\ref{frame}) and we still denote the corresponding lifting $\widetilde{\gamma}\in PHS(n)$ of $\gamma$ by
$$\widetilde{\gamma}(s)=(\gamma(s),e_1(s),\cdots,e_n(s),e_{n+1}(s),\cdots,e_{2n}(s),T), $$which is unique up to a $SO(2n)$ group action.
Let $\omega$ be the Maurer-Cartan form of $PSH(n)$. We shall derive the Darboux derivative $\widetilde{\gamma}^*\omega$ of the $\widetilde{\gamma}.$

By the moving frame formulas $(\ref{mformula})$,
$$d\gamma(s)=e_j(s)\widetilde{\gamma}^*\omega^j+e_{n+j}(s)\widetilde{\gamma}^*\omega^{n+j}+T\widetilde{\gamma}^*\omega^{2n+1};$$
on the other hand, by the choice of the oriented frame
$$d\gamma(s)=\gamma_{\xi}^{'}(s)ds+\gamma_T^{'}(s)ds=e_1(s)ds+\gamma_T^{'}(s)ds.$$
Comparing the components in above equations, we have
\begin{align}\label{dr}
\begin{array}{ll}
&\widetilde{\gamma}^*\omega^1=ds, \ \ \ \widetilde{\gamma}^*\omega^\ell=0\  \text{ for } 2\leq \ell \leq 2n,\\
& \widetilde{\gamma}^*\omega^{2n+1}
=\langle\frac{d\gamma(s)}{ds},T\rangle ds=\tau(s)ds.\end{array}
\end{align}
Again from (\ref{mformula}), we have
\begin{align}\label{dej1}
\begin{array}{ll}
de_j(s)&=e_k(s)\widetilde{\gamma}^*\omega^k_j+e_{n+k}(s)\widetilde{\gamma}^*\omega^{n+k}_j+T\widetilde{\gamma}^*\omega^{n+j}
\\& =e_k(s)\widetilde{\gamma}^*\omega^k_j+e_{n+k}(s)\widetilde{\gamma}^*\omega^{n+k}_j \end{array}
\end{align}
For $1\leq j\leq n-1$, since $e_j(s)$ (resp.  $de_j(s)$) is the linear combination of
$\gamma_{\xi}^{'}(s),\cdots,\gamma_{\xi}^{(j)}(s),$ (resp. $\gamma_{\xi}^{'}(s),\cdots,\gamma_{\xi}^{(j+1)}(s)$ ), one has
\begin{align}\label{dej2}
\widetilde{\gamma}^*\omega^i_j=0,\mbox{ for }\ \
i>j+1.
\end{align}
By (\ref{dej1}) and the definition of p-curvatures, we also have
\begin{align}\label{dej3}
\widetilde{\gamma}^*\omega_j^{j+1}=\langle\frac{de_j(s)}{ds},e_{j+1}(s)\rangle ds=\kappa_j(s)ds.
\end{align}
Similarly, for $j=n$,
$$\begin{array}{ll}de_n(s)&=e_k(s)\widetilde{\gamma}^*\omega^k_n+e_{n+k}(s)\widetilde{\gamma}^*\omega^{n+k}_n+T\widetilde{\gamma}^*\omega^{2n}
\\ &=e_k(s)\widetilde{\gamma}^*\omega^k_n+e_{n+k}(s)\widetilde{\gamma}^*\omega^{n+k}_n.\end{array}$$
In addiction, since $\omega^{n+k}_n=-\omega^{k}_{2n}=\omega_k^{2n}$,
$\widetilde{\gamma}^*\omega^{n+i}_n=\widetilde{\gamma}^*\omega_i^{2n}=0$, for $1\leq j\leq
n-1$. One has
\begin{align}\label{den}
\widetilde{\gamma}^*\omega_n^{2n}=\langle\frac{de_n(s)}{ds},e_{2n}(s)\rangle ds=\kappa_n(s)ds.
\end{align}

By (\ref{dr})(\ref{dej2})(\ref{dej3})(\ref{den}) and use the anti-symmetric property, $\omega^i_j=-\omega^j_i$, finally we reach the moving frame formulaes for the curve $\gamma(s)$
\begin{align}\label{mformular}
\begin{array}{rcccl}
d\gamma(s)&=&e_1(s)ds&+&T\tau(s)ds, \\
de_j(s) &=&-e_{j-1}(s)\kappa_{j-1}(s)ds&+&e_{j+1}(s)\kappa_j(s)ds, \\
de_n(s)&=&-e_{n-1}(s)\kappa_{n-1}(s)ds &+&e_{2n}(s)\kappa_n(s)ds.
\end{array}
\end{align}
In conclusion, we obtain the Darboux derivative $\widetilde{\gamma}^*\omega$ of $\widetilde{\gamma}$
\begin{align}\label{darboux}
&\widetilde{\gamma}^*\omega= \nonumber \\
&\left(\begin{array}{c|ccccc|ccccc|c}
0&0&\cdots&\cdots &\cdots &0&0&\cdots &\cdots&\cdots &0&0\\ \hline
1&0&-\kappa_1(s)& 0 &\cdots& 0 & 0 &\cdots  & \cdots & \cdots &0&0\\
0&\kappa_1(s)&\ddots& \ddots&\ddots&\vdots &\vdots &\ddots & &&\vdots&\vdots\\
\vdots&0&\ddots& \ddots&\ddots& 0&\vdots& &\ddots& &\vdots&\vdots\\
\vdots&\vdots&\ddots& \ddots&\ddots &-\kappa_{n-1}(s)&\vdots&  & &0 &0&\vdots\\
0&0&\cdots& 0 &\kappa_{n-1}(s)&0&0&\cdots&\cdots&0&-\kappa_n(s)&0\\ \hline
0&0&\cdots&\cdots&\cdots&0&0&-\kappa_1(s)& \cdots&0&0& 0\\
\vdots&\vdots&\ddots&& &\vdots&\kappa_1(s)&0&\ddots& \ddots&0&\vdots\\
\vdots&\vdots& &\ddots& &\vdots&\vdots&\ddots& \ddots&\ddots&\vdots&\vdots\\
\vdots&\vdots&\cdots&&0&0&0&\cdots& \ddots&0&-\kappa_{n-1}(s)&\vdots\\
0&0&\cdots&&0&\kappa_n(s)&0&\cdots& \cdots&\kappa_{n-1}(s)&0&0\\ \hline
\tau(s)&0&\cdots& \cdots&\cdots&0&0 &\cdots&\cdots &0&-1&0\\
\end{array}\right)ds.
\end{align}

\section{proof of the main theorem}
We show the proof of Theorem \ref{mthm1} in this section.
\begin{proof}
First we show the existence. Given $(n+1)$ functions $\kappa_i(s)$, $1\leq i\leq n,$ and $ \tau(s)$
defined on an open interval $I$. Define a $PSH(n)$-valued one-form $\varphi$ on $I$ with entries $k_i$'s, $\tau$ as the one in (\ref{darboux}).
It is easy to show that $d\varphi+\varphi\wedge \varphi=0$, and by Theorem \ref{existence}, there exists a curve
$$\widetilde{\gamma}(s)=(\gamma(s); e_1(s),\cdots,e_n(s),e_{n+1}(s),\cdots,e_{2n}(s),T)\in PSH(n)$$
such that $\widetilde{\gamma}^*\omega=\varphi$. Therefore, by the moving frame formula (\ref{mformular}), we have
\begin{align*}
d\gamma(s)&=e_1(s)ds+T\tau(s)ds, \\
de_j(s) &=-e_{j-1}(s)\kappa_{j-1}(s)ds+e_{j+1}(s)\kappa_j(s)ds, \\
de_n(s) &=-e_{n-1}(s)\kappa_{n-1}(s)ds +e_{2n}(s)\kappa_n(s)ds,\\
de_{n+j}(s)&=-e_{n+j-1}(s)\kappa_{j-1}(s)ds +e_{n+j+1}(s)\kappa_j(s)ds, \\
de_{2n}(s)&=-e_{n}(s)\kappa_{n}(s)ds -e_{2n-1}(s)\kappa_{n-1}(s)ds,
\end{align*} for $1\leq j \leq n-1$,
which implies that
\begin{align*}
e_1(s)&=\gamma_{\xi}^{'}(s),\\
\kappa_j(s)&=\langle\frac{de_j(s)}{ds},e_{j+1}(s)\rangle, \ 1\leq j\leq n-1,  \\
\kappa_n(s)&=\langle\frac{de_n(s)}{ds},e_{2n}(s)\rangle,\\
\tau(s)&=\langle\frac{d\gamma(s)}{ds}, T\rangle.
\end{align*}
We have reached the proof of existence.

Next, for the uniqueness, suppose that $\gamma_1$ and $\gamma_2$ have the same
$p$-curvatures $\kappa_j(s)$, $1\leq j\leq n$ and the contact normality $\tau(s)$.
By the moving frame formulas (\ref{mformula}) we get
$$\widetilde{\gamma}_1^*\omega=\widetilde{\gamma}_2^*\omega.$$
Therefore, Theorem \ref{uniqueness} implies that there exists an element $g \in PSH(n)$ such that
$\widetilde{\gamma}_2(s)=g\circ\widetilde{\gamma}_1(s)$, and hence
$\gamma_2(s)=g \circ \gamma_1(s)$ for all $s$. This completes the
proof of uniqueness.
\end{proof}

\section{The degenerate case}
We give the proof of Proposition \ref{topdeg}.
\begin{proof}
Without lose of generality, we may assume order($\gamma$)=$n-1$ and $\gamma(0)=(\beta(0),z(0))=0$. We observe that the second condition in $(\ref{order})$ holds if and only if $$\gamma'_\xi(s)\wedge \cdots \gamma_\xi^{(n-1)}(s)\wedge J\gamma_\xi(s)\wedge \cdots J\gamma_\xi^{(n-1)}(s)\neq 0.$$
At $s=0$, we may take the orthonormal frame $e_1(0), \cdots , e_{n-1}(0)$ satisfying
$e_1(0)\wedge \cdots \wedge e_k(0)=
\frac{\gamma_\xi'\wedge \cdots \wedge \gamma_\xi^{(k-1)}}{|\gamma_\xi'\wedge \cdots \wedge \gamma_\xi^{(k-1)}|}$ for all $1\leq k\leq n-1$
such that $span_\mathbb{R}\{e_1(0)$, $\cdots$, $e_{n-1}(0)$, $Je_1(0)$, $\cdots$, $Je_{n-1}(0)\}=\mathbb{C}^{n-1}\subset \mathbb{C}^n$. We also have the natural orthogonal decomposition $$\mathbb{C}^n\ni \beta(s)=\widetilde{\beta}(s)+a(s)N$$ for some function $a(s)$, where $N$ is normal to $\mathbb{C}^{n-1}$ and $\widetilde{\beta}\in \mathbb{C}^{n-1}$. Since $\beta(0)\in \mathbb{C}^{n-1}$, we have the initial condition
\begin{align}\label{initial}a(0)=0.\end{align}

We shall claim that $a(s)=0$ for all $s$, which implies that $\gamma(s)\in \mathbb{C}^{n-1}$.
Since $\beta^{(j)}(0)\in \mathbb{C}^{n-1}$,
\begin{align}\label{eqns}
a^{(j)}(0)=0 \text{ for }1\leq j\leq n-1.
\end{align}
On the other hand, by the assumption
\begin{align*}
0&=\beta'\wedge \cdots \wedge \beta^{(n)}\\
&=(\widetilde{\beta}'+a'N)\wedge \cdots \wedge (\widetilde{\beta}^{(n)}+a^{(n)}N)\\
&=(\widetilde{\beta}'\wedge \cdots \wedge \widetilde{\beta}^{(n)})+\sum_{h=1}^n  \widetilde{\beta}'\wedge\cdots\wedge \widetilde{\beta}^{(h-1)}\wedge (a^{(h)}N)\wedge \widetilde{\beta}^{(h+1)}\wedge \cdots \widetilde{\beta}^{(n)}\\
&=(a'b_1+a''b_2+\cdots+a^{(n-1)}b_{n-1}+a^{(n)})(\widetilde{\beta}'\wedge \cdots \widetilde{\beta}^{(n-1)}\wedge N),
\end{align*}where we denote $\widetilde{\beta}^{(n)}=\sum_{j=1}^{n-1}b_j\widetilde{\beta}^{(j)}$ for some smooth functions $b_j(s)$ independent of $a^{(j)}$'s and we have used the fact that the complex $n$-form $\bigwedge_{h=1}^n\tilde{\beta}^{(h)}$ vanishes in $\mathbb{C}^n$. Since the volume form $\big(\bigwedge_{h=1}^{n-1}\widetilde{\beta}^{(h)}\big)\wedge N$ is nonzero, together with (\ref{initial}), (\ref{eqns}), we obtain the $n$-th order O.D.E. system with the initial conditions
\begin{align*}
\left\{
\begin{array}{rl}
a^{(n)}(s)+\sum_{j=1}^{n-1}b_j(s)a^{(j)}(s)&=0, \\
a^{(j)}(0)&=0 \text{ for }j=0,\cdots, n-1.
\end{array}\right.
\end{align*}Hence by the existence and uniqueness theorem of O.D.E. one has $a(s)\equiv 0$ for all $s$, which implies $\gamma(s)\in H_{n-1}\subset H_n$, and complete the proof.
\end{proof}

Finally we calculate $order(\gamma)=1$ if $\gamma$ is a horizontal geodesic in $H_n$ for any $n\geq 1$.

\begin{proof}[\textbf{Proof of Proposition \ref{geodesic}}]
In \cite{Rit}, the horizontal geodesic $\gamma:I\rightarrow H_n$ satisfies the equation $$D_{\gamma'}
{\gamma'}+2\lambda J(\gamma')=0,$$ for some constant
$\lambda \in R.$

Let $\gamma(s)=(x_1(s),\cdots,x_n(s),y_1(s), \cdots,y_n(s), z(s))$ be a horizontal geodesic with horizontal arc-length $s$. Then, for $1\leq j\leq n$, we have the following expression
\begin{align}
x_j(s)&=(x_0)_j+A_j\Big(\frac{\sin(2\lambda s)}{2\lambda}\Big)+B_j\Big(\frac{1-\cos(2\lambda s)}{2\lambda}\Big), \label{xj} \\
y_j(s)&=(y_0)_j+A_j\Big(\frac{1-\cos(2\lambda s)}{2\lambda}\Big)+B_j\Big(\frac{\sin(2\lambda s)}{2\lambda}\Big), \label{yj}\\
z(s)&=t_0+\frac{1}{2\lambda}\Big(s-\frac{\sin(2\lambda s)}{2\lambda}\Big)+\sum^n_{j=1}\Big\{(A_j(x_0)_j+B_j(x_0)_j)(\frac{1-\cos(2\lambda
s)}{2\lambda}) \label{zs} \\
&\hspace{1cm}-(B_j(x_0)_j-A_j(x_0)_j)(\frac{\sin(2\lambda s)}{2\lambda})\Big\},  \nonumber
\end{align}
with the initial conditions $x_j(0)=(x_0)_j$, $y_j(0)=(y_0)_j$, and ${x'}_j(0)=A_j$, $y'_j(0)=B_j$ satisfying $\sum_{j=1}^n (A_j^2+B_j^2)=1$. By the decomposition (\ref{decomp})
\begin{align*}
\gamma^{'}(s)&=(x_1^{'}(s),\cdots,x_n^{'}(s),y_1^{'}(s),\cdots,y_n^{'}(s),z^{'}(s))\\
&=\sum_{j=1}^n \Big(x_j^{'}(s)\frac{\partial}{\partial x_j}+y_j^{'}(s)\frac{\partial}{\partial y_j}\Big)+z^{'}(s)\frac{\partial}{\partial
z}\\
&=\sum_{j=1}^n\Big( x_j^{'}(s)\mathring{e}_j+y_j^{'}(s)\mathring{e}_{n+j}\Big)+\sum_{j=1}^n \Big(z^{'}(s)+x_j y_j^{'}-y_jx_j^{'}\Big)\frac{\partial}{\partial
z},
\end{align*} we have
\begin{align}
&\gamma^{'}_{\xi}=\sum_{j=1}^n\Big(x_j^{'}(s)\mathring{e}_j+y_j^{'}(s)\mathring{e}_{n+j}\Big)
, \label{gamma'} \\ &
\gamma^{'}_{T}=\sum_{j=1}^n \Big(z^{'}(s)+x_iy_i^{'}-y_ix_i^{'}\Big)T, \nonumber
\end{align}
where $\frac{\partial}{\partial z}=T.$ Since the geodesic is horizontal, the contact normality $\tau(s)=0$. Moreover, by (\ref{xj}), (\ref{yj}), (\ref{gamma'})
\begin{align*}
\gamma^{'}_{\xi}(s)=\sum_{j=1}^n\Big( (A_j \cos(2\lambda s)+B_j \sin(2\lambda s))\mathring{e}_j+(-A_j \sin(2\lambda
s)+B_j \cos(2\lambda s))\mathring{e}_{n+j}\Big).
\end{align*}
Note that $|\gamma^{'}_{\xi}(s)|=1$, we may take $e_1=\gamma'_\xi$. Taking the derivatives, we observe that
\begin{align*}
&\gamma^{''}_{\xi}(s)=\sum_{j=1}^n \Big( (-2\lambda A_j \sin(2\lambda s)+2\lambda B_j \cos(2\lambda s))\mathring{e}_j+(-2\lambda A_j \cos(2\lambda
s)-2\lambda B_j \sin(2\lambda s))\mathring{e}_{n+j}\Big), \\
&\gamma^{'''}_{\xi}(s)=-(2\lambda)^2 \gamma^{'}_{\xi}(s),\\
&\langle \gamma^{''}_{\xi}(s), \gamma^{'}_{\xi}(s)\rangle =0, \mbox{ and } \gamma^{''}_{\xi}(s)=-2\lambda Je_1.
\end{align*}
By Definition \ref{deg},  we conclude that the order of geodesics is one. There is only one invariant, p-curvature, for $\gamma$, namely,
\begin{align*}
\kappa_1&=\langle \frac{de_1(s)}{ds}, e_2(s)\rangle \\
&=-2\lambda \sum_{j=1}^n \Big\{( A_j
\sin(2\lambda s)- B_j \cos(2\lambda s))^2+( A_j \cos(2\lambda s)+ B_j \sin(2\lambda s))^2 \Big\}\\
&=-2\lambda,
\end{align*}
where $e_2=Je_1$.
\end{proof}

\end{document}